\documentclass[12pt,twoside,a4paper]{amsart}
\usepackage[centertags]{amsmath}
\usepackage{amssymb,amscd,amsfonts,latexsym,a4wide,amsthm}
\usepackage{bbm}
\usepackage{times}
\newtheorem{thm}{Theorem}[section]
\newtheorem{cor}[thm]{Corollary}

\newtheorem{prop}[thm]{Proposition}

\theoremstyle{definition}

\newtheorem{prg}[thm]{}

\theoremstyle{remark}
\newtheorem{remark}[thm]{Remark}

\numberwithin{equation}{section}

\newcommand{\thmref}[1]{Theorem~\ref{#1}}

\newcommand{\propref}[1]{Proposition~\ref{#1}}

\newcommand{\corref}[1]{Corollary~\ref{#1}}

\newcommand{\formref}[1]{(\ref{#1})}
\newcommand{\prgref}[1]{{\bf\ref{#1}}}

\newcommand{\Cc}{{\mathbb C}}

\newcommand{\Zz}{{\mathbb Z}}
\newcommand{\Qq}{{\mathbb Q}}
\newcommand{\Pp}{{\mathbb P}}
\newcommand{\Rr}{{\mathbb R}}

\newcommand{\Spec}{\operatorname{Spec}}
\newcommand{\Ob}{{\mathcal O}}

\renewcommand{\leq}{\leqslant}
\renewcommand{\geq}{\geqslant}

\def\rep#1{\bysame}
\def\bib[#1]#2<#3>#4|#5(#6)#7-#8.{\bibitem{#1}
{\sc #2},\ #3,\ {\it #4},\ {\bf #5}\ (#6)\ #7--#8.}
\def\tbib[#1]#2<#3>#4.{\bibitem{#1} {\sc #2},\ #3,\ #4.}
\def\bibt[#1]#2<#3>#4|#5(#6)#7-#8{\bibitem{#1}
{\sc #2},\ #3,\ in\ {\it #4},\ ed.\ {#5},\ (#6)\ #7--#8.}
\def\bibit[#1]#2<#3>{\bibitem{#1} {\sc#2},\ #3}
\def\prebib[#1]#2<#3>{\bibitem{#1} {\sc #2},\ {\it #3}, Preprint.}
\def\toabib[#1]#2<#3>(#4){\bibitem{#1} {\sc #2},\ #3,\ Accepted for
publication in {\it #4}}

\title[Positivity of Heights]{Positivity of Heights of Semi-Stable Varieties}
\author[R.~G.~FERRETTI]{Roberto~G.~FERRETTI}
\address{Universit\'a della Svizzera italiana, Via Buffi 19, CH-6900 Lugano, Switzerland}
\address{Departement Mathematik, ETH Zentrum, CH-8092 Z\"urich, Switzerland}
\email{roberto.ferretti@lu.unisi.ch, ferretti@math.ethz.ch}
\date{\today}
\subjclass{11G50,14G40,14L24} \keywords{Arakelov Geometry, Geometric Invariant
Theory}
\begin{document}
\maketitle
\section{Introduction}
Let $K$ be a number field and ${\mathcal O}_K$ its ring of integers. Let
${\mathcal E}$ be an ${\mathcal O}_K$-module of rank $N+1$ in $\Pp({\mathcal
E}^\vee),$ the projective space representing lines in ${\mathcal E}^\vee.$ For
all closed subvarieties $X\subseteq \Pp({\mathcal E}_K^\vee)$ of dimension $d$
let $\deg(X)$ be its degree with respect to the canonical line bundle ${\mathcal
O}(1)$ of $\Pp({\mathcal E}_K^\vee).$ If ${\mathcal E}$ is endowed with the
structure of hermitian vector bundle over $\operatorname{Spec}({\mathcal O}_K)$
we can define the Arakelov degree $\widehat{\deg}\left(\overline{{\mathcal
E}}\right)$ and the Faltings height $h_{\mathcal E}(X)$ (see
\prgref{arakelov1}).\\
Let ${\mathcal V}$ be the direct limit of the set ${\mathcal V}_K$
of $\Ob_K$-modules ${\mathcal E}$ with an identity ${\mathcal
E}_K\cong K^{N+1}$ as $K$ varies in the set of number fields. Let
$X\subseteq \Pp^N$ be a closed irreducible projective variety of
dimension $d$ defined over $\overline{\Qq}.$ We define
$$
\hat{h}(X)=\inf_{{\mathcal E}\in {\mathcal
V}}\left(\frac{h_{\mathcal
E}(X)}{(d+1)\deg(X)}-\frac{\widehat{\operatorname{deg}}(\overline{{\mathcal
E}})}{N+1}\right).
$$
According to \cite[Th\'eor\`eme $1$]{Bo1} if the Chow point of $X$
is {\it semistable} (see \prgref{prg:chow}) then there exists a
non-negative constant $C$ such that
\begin{eqnarray}\label{bound}
\hat{h}(X)\geq -C.
\end{eqnarray}
This can be shown using the fact that for ${\mathcal
E}\in{\mathcal V}$ the expression $h_{\mathcal E}(X)/(d+1)\deg(X)$
is bounded below in terms of the height of a system of generators
of the rings of invariants
$\left(\operatorname{Sym}E\right)^{SL(E)}$ (see also \cite{Ga}).
Moreover, one can even show that $h_{\mathcal E}(X)/(d+1)\deg(X)$
is bounded below by the average of the successive minima of
${\mathcal E}$ (see \cite{So}). These results give a lower bounds
of $\hat{h}(X)$ that depends on the field $K$ of definition of
${\mathcal E}$ and $X.$ The independence of the constant $C$ from
the field of definition has been shown by Zhang in \cite{Z1},
using his absolute successive minima \cite[$5.2$]{Z3}, and
independently by Bost in \cite{Bo2}, using the ring of
invariants.\\
We will prove in this section a conjecture of Zhang
\cite[$4.3$]{Z1} that the height $\hat{h}(X)$ of a semistable
variety $X$ is non-negative, as it is true in the function field
case (see \cite{CH}, \cite{Bo1}). The idea of the proof is the
same of that for the main theorem of \cite{So}, the main
difference relies in the fact that here we make direct use of
Zhang's absolute successive minima of $X$ (defined in \cite{Z2})
at the place of Zhang's minima of the projective space.
\begin{thm}\label{th:princ}
Let us suppose that $X$ is semistable. Then
\begin{eqnarray}\label{bound2}
\hat{h}(X)\geq 0.
\end{eqnarray}
\end{thm}
In order to prove this inequality we have to introduce the {\it degree of
contact}. This is a birational invariant first considered in the context of
Geometric Invariant Theory (see \cite{Mu}), which has recently found important
applications in the domain of diophantine approximations. We refer to \cite{Mu},
and \cite[$\S 3$]{EF} for exhaustive and detailed discussions of the main
properties the degree of contact (see also \cite{FeCo},
\cite{FeDu} and the references therein).\\
\noindent As a byproduct, we are able to show, for all hermitian
${\mathcal O}_K$-module ${\mathcal E},$ new lower bounds for the
normalized Faltings height $h_{\mathcal E}(X),$ (see
\formref{normfal}). For instance, when $X$ is a generic $K3$
surface, i.e. a $K3$ surface whose Picard group has rank one, then
$X$ is semistable (see \corref{co:k3semistable} and \cite{Mo2}),
so \thmref{th:princ} holds true. This implies that the height
$h_{\mathcal E}(X)$ is non-negative for all hermitian ${\mathcal
O}_K$-module ${\mathcal E}.$ However, \corref{co:k3} provides a
lower bound for the Faltings height in terms of a linear
combination of Zhang's absolute minima which is stronger than that
given by the semistability of $X.$\\
The seminal work of Zhang \cite{Z1} has even contributed to lay
down the basis for the proof by Phong and Sturm \cite{PS} of a
conjecture formulated by Donaldson in his work \cite{Do} on Yau's
conjecture. However, Donaldson, Phong and Sturm were concerned
about one direction of Yau' conjecture, namely on the proof of the
semistability of $X$ under the assumption of the existence of a
K\"ahler-Einstein metric. It would be interesting to see if the
methodology developed in this note could be useful to better
understand the other direction of Yau's conjecture.
\section{Degree of contact}
\begin{prg}\label{prg:primo} Let $E$ be a finite dimensional vector space defined over a number field $K.$
A {\it weight function} is a map $w:E\to\Rr\cup\{-\infty\}$
satisfying:
\begin{enumerate}
\item $w(x)=-\infty,$ if and only if $x=0,$
\item for all $t\in K^*$ and all $x\in E,$ $w(t\cdot x)=w(x),$
\item for all $x,y\in E,$ $w(x+y)\leq \max\{w(x), w(y)\}$.
\end{enumerate} For all non-negative real numbers $\alpha$ the
set $ F^\alpha=\{x\in E\ :$ $\ w(x)\leq \alpha\}$ is a subspace of $E$, and
$F^\alpha\subseteq F^\beta$ whenever $\alpha\leq \beta.$ Varying $\alpha$, we
get then an exhaustive ($\cup_{\alpha\in \Rr}F^\alpha=E$) and separated
($\cap_{\alpha\in\Rr}F^\alpha=\{ 0\}$) filtration ${\mathcal F}$ of subspaces of
$E.$ We say that a basis $l_0,\cdots,l_N$ of $E$ is {\it adapted to the
filtration} ${\mathcal F}$, if for all $\alpha\in \Rr$
$$
F^\alpha=\bigoplus_{w(l_j)\leq \alpha}K\cdot l_j.
$$
We number the elements of this basis so that $w(l_0)\geq
w(l_1)\geq \cdots \geq w(l_N).$ For $h=0,\cdots,N$ define
$r_h:=w(l_h)$ and ${\bf r}=(r_0, \cdots,r_N).$ The vector ${\bf
r}$ is obviously independent of the chosen basis. This
construction gives a bijective map between the set of weight
functions on $E$ and the set of couples $({\mathcal F}, {\bf r}),$
with ${\mathcal F}$ a filtration of subspaces of $E,$ and ${\bf
r}=(r_0,$ $\cdots ,r_N)\in\Rr_{\geq 0}^{N+1}$ with
$r_0\geq\cdots\geq r_N$. Indeed, to such a pair $({\mathcal F},
{\bf r})$ we associate a weight function $w$ on $E$ as follows:
Put $w(0)=-\infty$ and for $x\in E\setminus\{0\}$ define $w(x)$ as
the smallest $r_i$ such that $x\in F^{r_i}.$ Equivalently, given a
basis $l_0,\cdots,l_N$ of $E$ adapted to the filtration ${\mathcal
F}$, write $x=x_0l_0+\cdots x_Nl_N.$ Then $w(x)$ is the greatest
$r_i$ for which $x_i\neq 0$. Notice that the value $w(x)$ does not
change much if we dilate ${\bf r}$ or perturb it a little bit. We
can always find an integer valued weight function $\tilde{w}$
supported on the same filtration ${\mathcal F}$ of $w,$ such that
for some sufficiently small $\varepsilon>0$ and some positive
integer $m$ we have $mw\leq\tilde{w}\leq m(1+\varepsilon)w.$ This
enables us to reduce most of the computations to weight functions
with integer values.
\end{prg}
\begin{prg} \label{functorial} Weight functions satisfy several
functorial relations. Let $w$ be a weight function on $E$ with
non-negative integer weights ${\bf r}\in\Zz^{N+1}$ and associated
filtration ${\mathcal F}.$ Consider a subspace $F\subseteq E.$ The
restriction $w|_F$ of $w$ on $F$ defines a weight function on $F.$
Further, $w$ induces a weight function on the quotient $E/F,$
mapping $l$ to the minimum of the weights of the elements $x\in E$
with $\pi(x)=l,$ where $\pi :E\to E/F$ is the canonical
projection. An element $h\in E^\vee$ is a linear functional
$h:E\to K.$ If $h\neq 0$ we define the weight of $h$ as minus the
weight of the line
$E/\operatorname{ker}(h).$\\
Given two vector spaces $E_1$, and $E_2$ over $K,$ endowed with
weight functions $w_1, w_2$ respectively, we define a weight
function $w$ on $E_1\oplus E_2$ by
\begin{eqnarray}\label{wsum}
e_1\oplus e_2\mapsto\max\{w_1(e_1), w_2(e_2)\}.
\end{eqnarray}
Moreover, on the tensor product $E_1\otimes E_2$ we define a
weight function $w$ by $w(e_1\otimes e_2)= w_1(e_{1}) +
w_2(e_{2}).$ Let $m$ be a positive integer. The symmetric group of
order $m$ operates on the $m$-th tensor power $E^{\otimes m}$ by
permuting the factors. The $m$-th symmetric power of $E,$ denoted
by $\operatorname{Sym}^mE,$ is the maximal subspace of $E^{\otimes
m}$ invariant under this operation. Whence by restriction we can
canonically define a weight function on $\operatorname{Sym}^mE.$
The exterior power $\bigwedge^m E$ is the quotient of $E^{\otimes
m}$ by $\operatorname{Sym}^mE,$ therefore on this space there
exists a canonical induced weight function.
\end{prg}
\begin{prg}\label{prg:hilbertchow} Let $w$ be a weight function on $E$ with integer
weights ${\bf r}\in\Zz^{N+1},$ and let ${\mathcal F}$ be the
associated filtration. Let $X$ be a projective absolutely
irreducible scheme of dimension $d$ and defined over $K$ embedded
into the projective space $\Pp(E^\vee)$ representing lines of the
dual vector space $E^\vee.$ If $m$ is large enough, say $m\geq
m_0$, the cup product map
\begin{eqnarray}\label{hilb}
\varphi_m: \operatorname{Sym}^m(E)\to H^0(X,{\mathcal O}(m))
\end{eqnarray}
is surjective. Therefore $H^0(X,{\mathcal O}(m))$ can be
identified with a quotient of $\operatorname{Sym}^m(E).$ As in
\prgref{functorial} $w$ induces then a weight function on
$H^0(X,{\mathcal O}(m)),$ and on the one-dimensional space
$\bigwedge^{h^0(X,\Ob(m))}H^0(X,{\mathcal O}(m)).$ We denote by
$w(X,m)$ the weight of this line, which is well defined by the
homogeneity property $(2)$ of weight functions. There exists an
integer $e_{w}(X)$ such that when $m$ goes to infinity
\begin{eqnarray}\label{asym}
w(X,m)& = & e_{w}(X)\frac{m^{d+1}}{(d+1)!}+O(m^d),
\end{eqnarray}
(see \cite[Proposition $2.11$]{Mu}, or \cite[$\S 3$]{EF}). The
number $e_{w}(X)$ is called {\it degree of contact} of $X$ with
respect to the weighted filtration associated to the weight
function $w$. We extend this definition by linearity to cycles,
and by approximating to real weights ${\bf r}\in \Rr^{N+1}.$\\
In the last years appeared several articles in diophantine
approximations that make a wide use of the degree of contact (see
\cite{EF},\cite{FeCo}, and \cite{FeDu}). In these articles the
main properties of the degree of contact are discussed in detail.
We refer to them for a further thorough analysis of the degree of
contact (see also \cite{Mu}, \cite{Mo1}, \cite{KSZ}).
\end{prg}
\begin{prg}\label{prg:chow}
Each suitably generic element ${\bf h}=(h_0,\cdots,
h_d)\in\Pp(E)^{d+1}$ defines naturally a $(N-d-1)$-dimensional
linear subspace $L_{\bf h}\subseteq \Pp(E^\vee).$ Consider the set
$Z(X)$ of all $(d+1)$-tuples ${\bf h}\in\Pp(E)^{d+1},$ such that
$L_{\bf h}(\overline{\Qq})\cap X(\overline{\Qq})\neq\emptyset,$
where $L_{\bf h}$ has dimension $N-d-1.$ Then $Z(X)$ is an
irreducible hypersurface of multidegree $(\deg(X),\cdots,\deg(X))$
(see \cite[Thm. IV, p. $41$]{HP}). The hypersurface $Z(X)$ turns
out to be given by an up to a constant unique polynomial element
$F_X\in V,$ where
\begin{eqnarray}
V & = &
\left[\left(\operatorname{Sym}^{\deg(X)}E\right)^{\otimes(d+1)}
\right]^\vee.\label{defv}
\end{eqnarray}
This is a so called {\it (Cayley-Bertini-van der Waerden-) Chow
form} of $X.$ By definition it has that property that
$F_X(h_0,\cdots,h_d)=0$ if and only if $X$ and the hyperplanes
given by the vanishing of the linear forms $h_i$ $(i=0,\cdots,d)$
have a point in common over $\overline{\Qq}.$
\end{prg}
\begin{prg} As in \prgref{functorial} $w$ induces a weight function on
$V$, again denoted by $w$. From \cite[Proposition $2.11$]{Mu} (see
also \cite[Theorem $4.1$]{EF}) we know that the degree of contact
corresponds to minus the weight of the Chow point:
\begin{eqnarray}\label{mum}
e_{w}(X)=-w(F_X).
\end{eqnarray}
Indeed, using the terminology of \cite{Mu}, the {\it ``n.l.c. of $r_n^V$,''}
(the degree of contact) corresponds to the {\it ``$\lambda$-weight $a_V$ of
$\phi_V$''} (minus the weight of the Chow form in our notation). \\
We say that the variety $X$ is {\it Chow-semistable} (or simply semistable) if
for all weight functions $w$ on $E$ we have
$$
\frac{e_{w}(X)}{(d+1)\deg(X)}\leq \frac{1}{N+1}\sum_{i=0}^Nr_i.
$$
According to the Hilbert-Mumford criterion (see \cite{Mu}), this is equivalent
to say that the Zariski closure of the orbit of a representative of $F_X$ in $V$
under $SL(E)$ does not contain $0.$
\end{prg}
\begin{prg}\label{pr:evertse} Let $l_{0},\cdots,l_N$ be a basis of
$E$ adapted to the filtration associated to the weight function
$w,$ which identifies then $E$ with $K^{N+1},$ and define
$T=\binom{N+1}{d+1}.$ Given the blocks of variables
$h_p=(h_{p0},\cdots,h_{pN})$ $(p=0,\cdots,d)$ we define for each
subset $I_k=\{ i_{k0},\cdots,i_{kd}\}$ of $\{0,\cdots,N\}$ with
$i_{k0}<\cdots< i_{kd}$ the {\it bracket} $[I_k]=[i_{k0}\cdots
i_{kd}]=\det(h_{p,i_{kq}})_{p,q=0,\cdots,d},$ for $k=1,\cdots,T.$
From \cite[Thm. IV, p.$41$]{HP} it follows that the Chow form
$F_X$ can be expressed as a polynomial in terms of such brackets.
We expand $F_X$ as a sum of monomials of brackets
$$
F_X=\sum_{{\bf j}=(j_1,\cdots,j_T)\in {\mathcal J}}a_{\bf
j}[I_1]^{j_1}\cdots [I_T]^{j_T}
$$
where $a_{\bf j}\neq 0,$ and $|{\bf j} |=\deg(X)$ for ${\bf j}\in
{\mathcal J}.$ Then if $w$ is a weight function given by the
weights $r_0\geq \cdots \geq r_N\geq 0$ we have
$$
e_{w}(X)=\min_{{\bf j}\in {\mathcal J}} \sum_{i=1}^Tj_i\left(
\sum_{k\in I_i}r_k\right).
$$
\end{prg}
\begin{prg}
Suppose that the hyperplanes defined by the vanishing of the
linear forms $l_{N-d},\cdots,l_{N}$ do not have a common point on
$X$ defined over $\overline{\Qq}.$ If $X$ is linear, this means
precisely that the restrictions to $X$ of the linear forms
$l_{N-d},\cdots,l_{N}$ are linearly independent. We notice that
for given real numbers $0\leq c_0 \leq \cdots\leq c_N$ and for
${\bf j}\in {\mathcal J}$ we have
\begin{eqnarray}\label{linindep}
\sum_{i=1}^Tj_i\left( \sum_{k\in I_i}c_k\right) & \leq &
\deg(X)\Big(c_{N-d}+\cdots+c_N\Big).
\end{eqnarray}
Further, due to our assumption that $X$ and the zero set of
$l_{N-d},\cdots,l_N$ do not meet, we have
$F_X(l_{N-d},\cdots,l_N)\neq 0.$ Therefore, $F_X$ must contain the
monomial $[N-d\cdots N]^{\deg(X)}.$ This implies that among the
terms $\sum_{i=1}^Tj_i\left( \sum_{k\in I_i}c_k\right),$ with
${\bf j}\in {\mathcal J}$ we have $\deg(X)(c_{N-d}+\cdots+c_N).$
But this is the largest among all the terms \formref{linindep},
whence $$ \max_{{\bf j}\in {\mathcal J}} \sum_{i=1}^Tj_i\left(
\sum_{k\in I_i}c_k\right)=\deg(X)\Big(c_{N-d}+\cdots+c_N\Big).
$$
Remember that $w$ is a weight function given by the basis
$l_0,\cdots,l_N$ and weights $r_0\geq\cdots\geq r_N\geq 0.$ We
define $c_i=r_0-r_i,$ for $i=0,\cdots,N.$ We have
\begin{eqnarray}
e_w(X) & = & \min_{{\bf j}\in {\mathcal J}} \sum_{i=1}^Tj_i\left( \sum_{k\in I_i}r_k\right)\nonumber \\
& = & \deg(X)(d+1)r_0-\max_{{\bf j}\in {\mathcal J}} \sum_{i=1}^Tj_i\left( \sum_{k\in I_i}c_k\right)\nonumber \\
& = & \deg(X)(d+1)r_0-\deg(X)\Big(r_0-r_{N-d}+\cdots+r_0-r_{N}\Big)\nonumber \\
& = & \deg(X)\Big(r_{N-d}+\cdots+r_N\Big).\label{contlinind}
\end{eqnarray}
\end{prg}
\begin{prg} Let $X$ be an absolute irreducible projective surface and $C$ a
{\it pseudo ample} divisor on $X$, i.e. a divisor such that the
linear series $|C|$ has no fixed components and the associated map
$\phi_C:X\to\Pp(H^0(X,mC)^\vee)$ is birational onto its image, for
$m$ sufficiently large. For $E=H^0(X,mC)$ with $\dim(E)=N+1$ let
$w$ be a weight function on $E$ whith associated filtration
\begin{eqnarray}\label{base} {\mathcal F}: & E=V_0\supseteq
V_1\supseteq\cdots\supseteq V_{N}\supseteq\{0\}.
\end{eqnarray}
Suppose that there exists a blow up $\pi:B\to X$ on which $C$ has
a proper transform $\tilde{C},$ and such that for each
$i=0,\cdots,N$ the pullbacks of the sections in $V_i$ generate an
invertible sub-sheaf $\Ob_B(C_i)$ of $\Ob(\tilde{C})$. According
to \cite[$4.4$]{Fu} $\Ob_B(C_i)=\pi^*\Ob(C)\otimes\pi^{-1}(J_i),$
where $J_i$ is the ideal sheaf defining the base locus of $|V_i|.$
This means that the number $C_i^2$ is the degree of the projection
of $X$ onto $\Pp(V_i^\vee).$ Define $$ e_{ij}=C^2-C_i\cdot
C_j,\quad e_j=C^2-C_j^2. $$ These numbers are independent of the
choice of $B.$ Note that for all $i=0,\cdots,N$ the number $e_i$
measures the drop in degree on projection to $\Pp(V_i^\vee).$ Let
us identify $X$ with the image of the birational map $\phi_C$. Let
$r_0\geq\cdots\geq r_N\geq 0$ be the weights associated to the
weight function $w$. Then for all sequences of integers
$J=(j_0,\cdots,j_l)$ with $0=j_0<$ $j_1<\cdots<$ $j_l=N $
\cite[Proposition $2.10$]{FeDu} yields
\begin{eqnarray}\label{sj}
e_{w}(X) \leq
\sum_{k=0}^{l-1}(r_{j_k}-r_{j_{k+1}})\Big(e_{j_k}+e_{j_k,
j_{k+1}}+e_{j_{k+1}}\Big)=:S_J.
\end{eqnarray}
Assume now that $X$ is a $K3$ surface over $K,$ $C$ a pseudo ample
divisor on $X,$ and $m=1.$ Then
\begin{eqnarray}
h^0(X,C)=\frac{C^2}{2}+2\quad\text{ and }\quad
h^1(X,C)=0.\label{cinque}
\end{eqnarray}
Further $C$ has no  base points and $\phi_C(X)$ is projectively
normal (\cite[$2.6$]{SD}, \cite[$3.2$]{SD}, \cite[$6.1$]{SD}). If
we assume that no curve is contained in the base locus of $V_i$,
for all $i=0,\cdots,N,$ and that $r_N=0,$ then from \cite[Lemma
$5$]{Mo2} and \cite[Proposition $2.11$]{FeDu} we have
\begin{eqnarray}\label{degcontk3}
e_{w}(X) & \leq & -4r_0+6\sum_{i=0}^{N}r_i.
\end{eqnarray}
We are moreover able to prove the following tightening of
\formref{degcontk3}.
\begin{prop}\label{pr:tightdegcontk3}
Assume now that $X$ is a $K3$ surface over $K,$ $C$ a pseudo ample
divisor on $X,$ and that no curve is contained in the base locus
of $V_i$, for all $i=0,\cdots,N.$ Then if $r_N=0$ we have
$$
e_w(X)\leq \min\left\{-4r_0+6\sum_{i=0}^{N}r_i,2(N-1)r_0\right\}.
$$
\end{prop}
\proof Remember that for $i=0,\cdots,N,$ the number $e_i$ is the
amount by which the degree of $X$ in $\Pp(E^\vee)$ is greater that
that of its image under the projection onto $\Pp(V_i^\vee).$ In
any case
\begin{eqnarray*}
e_i & \leq & C^2.
\end{eqnarray*}
By Riemann-Roch \formref{cinque} we have
$N+1=h^0(X,C)=\frac{C^2}{2}+2,$ which implies
\begin{eqnarray}\label{eic2}
e_i & \leq & 2(N-1).
\end{eqnarray}
Let us now consider the inequality \formref{sj} with $J=(0,N).$
Then from \formref{eic2} we get
$$
e_{w}(X) \leq S_J=(r_0-r_N)\Big(e_0+e_{0N}+e_N\Big)\leq 2(N-1)r_0
$$
Together with \cite[Proposition $2.11$]{FeDu} this concludes the
proof.\qed

We recover here the main result of \cite{Mo2}.
\begin{cor}\label{co:k3semistable} Let $X$ be a $K3$
surface whose Picard group has rank $1$ and $C$ be a primitive
divisor class on $X.$ Then $X$ is semistable.
\end{cor}
\proof Let $w$ be any weight function on $E=H^0(X,C)$ and $r_0\geq
\cdots\geq r_N\geq 0$ be the associated weights. We can assume
without restriction that $r_N=0.$ Let us first suppose that
$$
r_0\geq \frac{3}{N+1}\sum_{i=0}^{N}r_i.
$$
Together with Riemann-Roch's formula \formref{cinque} and
\propref{pr:tightdegcontk3} this implies
\begin{eqnarray*}
\frac{e_{w}(X)}{(\dim X +1)\deg(X)} & \leq &
\frac{-4r_0+6\sum_{i=0}^{N}r_i}{6(N-1)}\leq
\frac{1}{N+1}\sum_{i=0}^{N}r_i.
\end{eqnarray*}
We assume now
$$
r_0\leq \frac{3}{N+1}\sum_{i=0}^{N}r_i.
$$
Then again Riemann-Roch's formula \formref{cinque} and
\propref{pr:tightdegcontk3} imply
\begin{eqnarray*}
\frac{e_{w}(X)}{(\dim X +1)\deg(X)} & \leq &
\frac{2(N-1)r_0}{6(N-1)}\leq \frac{1}{N+1}\sum_{i=0}^{N}r_i,
\end{eqnarray*}
which concludes the proof.\qed
\begin{remark}\label{re:pic1} As remarked in \cite{Mo2},
since the generic member of the moduli space of $K3$ surfaces has
Picard group of rank $1$, this result covers almost all $K3$
surfaces.
\end{remark}
\end{prg}
\section{Arakelov Geometry}
\begin{prg}\label{arakelov1}
Let $K$ be a number field and let $\Ob_K$ be its ring of integers, and let
$S_\infty$ be the set of complex embeddings of $K$. If ${\mathcal M}$ is a
torsion-free $\Ob_K$-module of finite rank such that, for all $\sigma\in
S_\infty,$ the corresponding complex vector space $M_\sigma={\mathcal
M}\otimes_{\Ob_K}\Cc$ is equipped with a norm $|\cdot |_\sigma,$ we may think of
${\mathcal M}$ as a free $\Zz$-module equipped with the norm $|\cdot |$on
$M_\sigma={\mathcal M}\otimes_{\Ob_K} \Cc=\oplus_{\sigma\in S_\infty}M_\sigma$
defined by $|\sum_{\sigma\in S_\infty}x_\sigma|=\sup_{\sigma\in
S_\infty}|x_\sigma|_\sigma$ for $x_\sigma \in M_\sigma,$ $\sigma\in S_\infty.$
In particular, consider an hermitian vector bundle $\overline{E}=({\mathcal
E},h)$ over $\Spec(\Ob_K)$ in the sense of \cite{GS}. In other words, ${\mathcal
E}$ is a torsion-free $\Ob_K$-module of rank $N+1<\infty$, and for all
$\sigma\in S_\infty$, $E_\sigma$ is equipped with a hermitian scalar product
$h$, compatible with the isomorphism $E_\sigma\cong E_{\overline{\sigma}}$
induced by complex conjugation. We will then denote by $\|\cdot\|_\sigma$ the
norm on $E_\sigma$ and $\|\cdot\|$ the norm on ${\mathcal E}\otimes_\Zz\Cc$ as
above. If $N=0$ then the {\it Arakelov degree} of $\overline{{\mathcal E}}$ is
defined by
$$ \widehat{\deg}\left(\overline{{\mathcal E}}\right)=\log(\# ({\mathcal E}/s\cdot
{\mathcal
O}_K))-\sum_{\sigma\in S_\infty}\log\| s\|_\sigma, $$ where $s$ is any non zero
element of ${\mathcal E}$. In general, we define the {\it (normalized) Arakelov
degree} of $\overline{{\mathcal E}}$ as $\widehat{\deg}_n(\overline{{\mathcal
E}})=\frac{1}{[K:\Qq]}\widehat{\deg}
\left(\det(\overline{{\mathcal E}})\right).$\\
Let ${\mathcal E}^\vee$ be the dual $\Ob_K-$module of ${\mathcal E},$ and let
$\Pp({\mathcal E}^\vee)$ be the associated projective space representing lines
in ${\mathcal E}^\vee.$ Consider a closed subvariety $X\subseteq \Pp(E^\vee),$
where $E={\mathcal E}\otimes K,$ of dimension $d,$ and let $\deg(X)$ be its
(algebraic) degree with respect to the canonical line bundle $\Ob(1)$ on
$\Pp(E^\vee)$. Let $h_{\mathcal E}(X)\in\Rr$ be the normalized Faltings height
of the Zariski closure $\overline{X}$ of $X$ in $\Pp({\mathcal E}),$ denoted by
$h_F(X)/[K:\Qq]$ in \cite[$(3.1.1),$ $(3.1.5)$]{BGS}. Let $\overline{\Ob(1)}$ be
the canonical line bundle equipped with the metric induced by $h$, then
\begin{eqnarray}\label{normfal}
h_{\mathcal E}(X) & = &
\frac{1}{[K:\Qq]}\widehat{\operatorname{deg}}\left(\hat{c}_1(
\overline{\Ob(1)})^{d+1}|\overline{X}\right)\in\Rr,
\end{eqnarray}
where $(.|.)$ is the bilinear pairing defined in loc. cit.
\end{prg}
\begin{prg} Let $L$ be a finite field extension of $K,$ and let $\overline{\mathcal
L}_0,\cdots,\overline{\mathcal L}_N$ be hermitian line subbundles of
$\overline{\mathcal E}$ such that $(\oplus_i{\mathcal L}_i)_L$ generates
${\mathcal E}_L.$ There exist points $P_0\cdots,P_N\in \Pp({\mathcal E}^\vee_L)$
associated to the line bundles above such that
$h(P_i)=-\widehat{\operatorname{deg}}(\overline{\mathcal L}_i)$ for
$i=0,\cdots,N,$ \cite[Theorem $5.2$]{Z2}. This implies that all
$\widehat{\operatorname{deg}}(\overline{\mathcal L}_i)$ are non-positive. Assume
that these line bundles are ordered by increasing Arakelov degree $
\widehat{\operatorname{deg}}(\overline{\mathcal L}_0)\leq\cdots$
$\leq\widehat{\operatorname{deg}}(\overline{\mathcal L}_N)\leq 0.$ For
$i=0,\cdots,N$ define $s_i=-\widehat{\operatorname{deg}}(\overline{\mathcal
L}_{i})+\widehat{\operatorname{deg}}(\overline{\mathcal L}_N),$ and put ${\bf
s}=(s_0,\cdots,s_N)\in\Rr^{N+1}.$\\
Let $x_0,\cdots,x_N$ be nonzero sections of the line bundles
$\overline{{\mathcal L}_0},\cdots,$ $\overline{{\mathcal L}_N},$ respectively,
that give an embedding $X\to\Pp^N.$\\ Further, let $w$ be a weight function on
$E={\mathcal E}_K$ with weights $r_0\geq \cdots\geq r_N=0,$ and ${\bf
r}=(r_0,\cdots,r_N).$ We get the following variation of \cite[Theorem $1$]{So}:
\begin{thm}\label{th:uno}
Assume there exists a continuous function $\psi:\Rr^{N+1}\to\Rr$ such that
$\psi(tx)=t\psi(x)$ for all $t\in\Rr,$ $x\in\Rr^{N+1},$ and such that $$
e_{w}(X) \leq \psi({\bf r}). $$ Then the following inequality holds: $$
h_{\mathcal E}(X)\geq
(d+1)\deg(X)\widehat{\operatorname{deg}}(\overline{{\mathcal L}}_N)-\psi({\bf
s}). $$
\end{thm}
{\it Proof: } This is a reinterpretation of the proof of  of \cite[Theorem
$4.4$]{Z1} (see \cite[$4.8$]{Z1}) using the language of  of \cite[Theorem
$1$]{So}.\qed
\begin{cor}\label{co:mixed} Let $X$ be a surface and ${\mathcal E}$ be a ${\mathcal
O}_K$-module with ${\mathcal E}_K=H^0(X,\Ob(C)),$ where $C$ is a pseudo ample
divisor. Then
$$
h_{\mathcal E}(X)\geq
(d+1)\deg(X)\widehat{\operatorname{deg}}(\overline{{\mathcal
L}}_N)-\sum_{k=0}^m\left( \widehat{\operatorname{deg}}(
\overline{{\mathcal
L}}_{j_{k+1}})-\widehat{\operatorname{deg}}(\overline{{\mathcal
L}}_{j_{k}})\right)\Big(e_{j_k}+e_{j_kj_{k+1}}+e_{j_{k+1}}\Big).
$$
\end{cor}
\proof Straightforward consequence of \thmref{th:uno} and
\formref{sj}.\qed
\begin{cor}\label{co:k3} Assume that $X$ is a $K3$ surface, and ${\mathcal E}$ be a
${\mathcal O}_K$-module with ${\mathcal E}_K=H^0(X,\Ob(C)),$ where $C$ is a
pseudo ample divisor. For $i=0,\cdots,N,$ let $V_i$ be the linear space of the
filtration \formref{base}. Suppose that no curve is contained in the base locus
of $V_i,$ for all $i=0,\cdots,N.$ Then
\begin{eqnarray}
\frac{1}{2}h_{\mathcal E}(X) & \geq &
\max\Big\{(N-1)\Big(\widehat{\operatorname{deg}}(\overline{{\mathcal
L}}_0)+2\widehat{\operatorname{deg}}(\overline{{\mathcal
L}}_N)\Big),\nonumber\\
&& -2\Big(\widehat{\operatorname{deg}}(\overline{{\mathcal
L}}_0)+2\widehat{\operatorname{deg}}(\overline{{\mathcal
L}}_N)\Big)+3\sum_{i=0}^{N}\widehat{\operatorname{deg}}(\overline{{\mathcal
L}}_i)\Big\}. \label{heightk3}
\end{eqnarray}
\end{cor}
{\it Proof: } From the assumption that  no curve is contained in
the base locus of $V_i$, for all $i=0,\cdots,N,$ given in
\formref{base} we know that we can estimate the degree of contact
with the help of the estimate \formref{pr:tightdegcontk3}. By
\thmref{th:uno} and Riemann-Roch \formref{cinque} we get
$$
h_{\mathcal E}(X)\geq
6(N-1)\widehat{\operatorname{deg}}(\overline{{\mathcal L}}_N)-\min
\left\{ -4s_0+6\sum_{i=0}^{N}s_i,2(N-1)s_0\right\},
$$
where $s_i=\widehat{\operatorname{deg}}(\overline{{\mathcal
L}}_N)-\widehat{\operatorname{deg}}(\overline{{\mathcal L}}_i).$
Expanding of this formula we get \formref{heightk3}, which
concludes the proof.\qed
\begin{cor}\label{co:notmeet}
Assume that the linear space defined by the vanishing of the last $d+1$ sections
$x_{N-d},\cdots,x_N$ does not meet $X.$ Then
$$
h_{\mathcal E}(X)\geq
\deg(X)\sum_{i=N-d}^N\widehat{\operatorname{deg}}(\overline{{\mathcal L}}_i)
$$
\end{cor}
\proof Follows from the identity \formref{contlinind} and \thmref{th:uno}.\qed
\end{prg}
\begin{prg}
We recall the definition of the normalized height $\hat{h}(X)$
form the introduction:
$$
\hat{h}(X)=\inf_{{\mathcal E}\in {\mathcal V}}\left(\frac{h_{\mathcal
E}(X)}{(d+1)\deg(X)}-\frac{\widehat{\operatorname{deg}}(\overline{{\mathcal
E}})}{N+1}\right),
$$
where ${\mathcal V}$ is the direct limit of the set ${\mathcal
V}_K$ of $\Ob_K$-modules ${\mathcal E}$ with an identity
${\mathcal E}_K\cong K^{N+1}$ as $K$ varies in the set of number
fields. We will prove in this section that under some conditions
on $X$ the height $\hat{h}(X)$ is non-negative, as it is true in
the function field case (see \cite{CH}, \cite{Bo1}). The idea of
the proof is the same of \cite{So}, the main difference relies in
the direct use of the minima of $X$ at the place of the minima of
the projective space.
\end{prg}
\begin{prg} According to \cite[$(5.2)$]{Z1} (see also \cite[Th\'eor\`eme $3.1$]{DP})
we know that for all $\varepsilon>0$ the set
\begin{eqnarray}\label{le:chowfalt}
X\left(\frac{h_{\mathcal E}(X)}{\deg(X)}+\varepsilon\right) & = & \left\{x\in
X\left(\overline{\Qq}\right)\Big|\ h(x)\leq \frac{h_{\mathcal
E}(X)}{\deg(X)}+\varepsilon \right\}
\end{eqnarray}
is Zariski dense.
\begin{prop}\label{pr:quasifine}
Assume that $X$ is not contained in any proper linear subspace of $\Pp(E^\vee).$
Then for all $\varepsilon>0$ there exist a finite field extension $L$ of $K$ and
$N+1$ hermitian line subbundles $\overline{{\mathcal
L}}_0,\cdots,\overline{{\mathcal L}}_N$ of $\overline{{\mathcal E}}$ such that
\begin{enumerate}
\item $\oplus_{i=0}^N{\mathcal L}_i$ generates generically ${\mathcal E}$ overl
$L,$
\item The arithmetic degrees of this subbundles satisfy
\begin{eqnarray}\label{chowempty2}
-\frac{1}{N+1}\sum_{i=0}^{N}\widehat{\operatorname{deg}}(\overline{{\mathcal
L}}_i)\leq\frac{h_{\mathcal E}(X)}{\deg(X)}+\varepsilon.
\end{eqnarray}
\end{enumerate}
\end{prop}
\proof Consider the $(N+1)$-fold Segre immersion $\varphi$ of $X^{N+1}\subseteq
\Pp(E^\vee)^{N+1}$ into the projective space $\Pp((E^\vee)^{\otimes (N+1)}).$
From \cite[$(2.3.19)$]{BGS} we have
$$
h_{(\mathcal E)^{\otimes
(N+1)}}\left(\varphi_*(X^{N+1})\right)=(N+1)\deg(X)^Nh_{\mathcal E}(X),
$$
and
$$
\deg\left(\varphi_*(X^{N+1})\right)=\deg(X)^{N+1}.
$$
Applying \formref{le:chowfalt} to $\varphi_*(X^{N+1})$ we have that the set
$\varphi_*(X^{N+1})\left(\frac{h_{\mathcal
E}(\varphi_*(X))}{\deg(\varphi_*(X))}+\varepsilon\right)$ is Zariski dense.
Since $\varphi$ is an isomorphism, this set is in homeomorphic to
$$
X^{N+1}\left((N+1)\frac{h_{\mathcal E}(X)}{\deg(X)}+\varepsilon\right),
$$
which is therefore Zariski dense, too. Further, notice that the
condition $(1)$ is obviously open, since it corresponds to the
non-vanishing of the maximal exterior power of the direct sum of
the line subbundles. Hence $(1)$ defines a Zariski open subset of
$$X^{N+1}\left((N+1)\frac{h_{\mathcal E}(X)}{\deg(X)}+\varepsilon\right),$$ which
is itself Zariski dense. This implies that the set of points with $(1)$ and
$(3)$ is Zariski dense, hence non-empty. Let $P=(P_0,\cdots,P_N)$ be a point
there, and $L$ be the field $ K(P_0,\cdots,P_N),$ then $P\in
\Pp(E^\vee)^{N+1}(L).$ Further, we let ${\mathcal L}_i$ the line subbundle of
${\mathcal E}$ associated to $P_i,$ $i=0,\cdots,N.$ By construction they satisfy
$(1)$ and $(2).$ This concludes the proof of
the Proposition.\qed\\
\noindent {\it Proof\ of \thmref{th:princ}.} According to \cite[Proposition
$4.2$]{Z1} it suffices to prove that, for any number field $K,$ $h_{\mathcal
E}(X)/((d+1)\deg X)\geq 0$ for any hermitian ${\mathcal O}_K$-module
$\overline{\mathcal E}$ such that $\det{\mathcal E}=0.$ By
\propref{pr:quasifine} for any $\varepsilon>0$ there exist a finite field
extension $L$ of $K$ and $N+1$ hermitian line subbundles $\overline{\mathcal
L}_0,\cdots,\overline{\mathcal L}_N$ of $\overline{\mathcal E}$ such that
$\oplus_{i=0}^N{\mathcal L}_i$ generates ${\mathcal E}$ over $L$ and
\formref{chowempty2} holds. For each $i=0,\cdots,N$ let $x_i$ be a section of
${\mathcal L}_i.$ These sections give an embedding $X\to\Pp^N.$ Assume that the
line bundles are ordered by increasing arithmetic degree
$$
\widehat{\operatorname{deg}}(\overline{{\mathcal L}}_0)\leq \cdots \leq
\widehat{\operatorname{deg}}(\overline{{\mathcal L}}_N)\leq 0.
$$
Since $X$ is semistable with respect to the sections
$x_0,\cdots,x_N$ and to any $(N+1)$-tuple of integers with
$r_0\geq \cdots\geq r_N\geq 0,$ we have
$$
\frac{e_w(X)}{(d+1)\deg(X)}\leq\frac{1}{N+1}\sum_{i=0}^Nr_i.
$$
Hence, from \thmref{th:uno} we get
\begin{eqnarray*}
\frac{h_{\mathcal E}(X)}{(d+1)\deg(X)} & \geq &
\widehat{\operatorname{deg}}(\overline{\mathcal
L}_N)+\frac{1}{N+1}\sum_{i=0}^N\Big(\widehat{\operatorname{deg}}(\overline{\mathcal
L}_i)-\widehat{\operatorname{deg}}(\overline{\mathcal
L}_N)\Big)\\
& = & \frac{1}{N+1}\sum_{i=0}^N\widehat{\operatorname{deg}}(\overline{\mathcal
L}_i).
\end{eqnarray*}
Since $\overline{\mathcal L}_0,\cdots \overline{\mathcal L}_N$
satisfy \formref{chowempty2} this yields
$$
\frac{h_{\mathcal E}(X)}{(d+1)\deg(X)} \geq -\frac{h_{\mathcal
E}(X)}{\deg(X)}-\varepsilon,
$$
whence
$$
\frac{h_{\mathcal E}(X)}{(d+1)\deg(X)}\geq -\frac{\varepsilon}{d+2}.
$$
But $\varepsilon>0$ can be chosen arbitrarily small, hence this concludes the
proof of \formref{bound2}. \qed
\end{prg}
\noindent {\it Aknowledgements.} It is a great pleasure to thank
Jean-Beno\^it Bost for sharing with me his ideas on Arithmetic and
Geometric Invariant Theory. Thanks are also due to Christophe
Soul\'e for many useful comments and Jan-Hendrik Evertse for his
constant encouragement and the explanation of the proof of
\prgref{pr:evertse}.

\end{document}